\documentclass[12pt]{amsart}
\usepackage{amssymb}
\usepackage{pstricks}
\allowdisplaybreaks

\def\modd{{m\,\rm odd}}
\def\meven{{m\,\rm even}}

\def\WPi{{W_{\Pi}}}
\def\GPi{\Gamma_{\Pi}}
\def\CPi{C_{\Pi}}

\def\WI{W_{\!I}}
\def\WIp{W_{\!I'}}

\def\WJ{W_{\!J}}

\def\GJ{\Gamma_{\!J}}

\def\Wmin{W_{\rm \!min}}
\def\Cmin{C_{\rm min}}
\def\Cmino{\Cmin^\circ}
\def\Gmin{\Gamma_{\rm \!min}}

\def\WOmega{W_{\Omega}}
\def\COmega{C_{\Omega}}
\def\GOmega{\Gamma_{\Omega}}

\def\Wmax{W_{\rm \!max}}
\def\Cmax{C_{\rm max}}
\def\Cmaxo{\Cmax^\circ}
\def\Gmax{\Gamma_{\rm \!max}}

\def\sset{\subseteq}
\def\iso{\cong}
\def\aut{\mathop{\rm Aut}\nolimits}

\def\Jperp{J^\perp}
\def\Jpperp{{J'}^\perp}
\def\a{\alpha}
\def\aperp{\alpha^\perp}

\def\tildeE{\tilde{E}}
\def\Z{\mathbb{Z}}

\newtheorem{theorem}{Theorem}
\newtheorem{lemma}[theorem]{Lemma}
\theoremstyle{remark}
\newtheorem*{remark}{Remark}
\newtheorem*{remarks}{Remarks}

\def\Dscaleunit{.0035cm}
\def\Dscalenoderadius{20}
\def\Dscalesinglewidth{15}
\def\solid#1#2{\qdisk(#1,#2){\noderadius}}
\def\hollow#1#2{\pscircle[fillstyle=solid,fillcolor=white](#1,#2){\noderadius}}
\def\single#1#2#3#4{\psline[linewidth=\singlewidth](#1,#2)(#3,#4)}
\def\dashed#1#2#3#4{\psline[linewidth=\singlewidth,linestyle=dashed,dash=10 10](#1,#2)(#3,#4)}

\begin {document}
\title{Normalizers of Parabolic Subgroups of Coxeter Groups}
\author{Daniel Allcock}
\address{Department of Mathematics\\University of Texas, Austin}
\email{allcock@math.utexas.edu}
\urladdr{http://www.math.utexas.edu/\textasciitilde allcock}
%
\date{May 2, 2011}

\begin{abstract}
We improve a bound of Borcherds on the virtual cohomological dimension
of the non-reflection part of the normalizer of a parabolic subgroup
of a Coxeter group.  Our bound is in terms of the types of the
components of the corresponding Coxeter subdiagram rather than the
number of nodes.  A consequence is an extension of Brink's result that
the non-reflection part of a reflection centralizer is free.  Namely,
the non-reflection part of the normalizer of parabolic subgroup of
type $D_5$ or $A_\modd$ is either free or has a free subgroup of
index~$2$.
\end{abstract}

\maketitle

Suppose $\Pi$ is a Coxeter diagram, $J$ is a subdiagram and $\WJ\sset
\WPi$ is the corresponding inclusion of Coxeter groups.  The
normalizer $N_{\WPi}(\WJ)$ has been described in detail by Borcherds
\cite{Bo} and Brink-Howlett \cite{BH}.  
Such normalizers have significant applications to working out the
automorphism groups of Lorentzian lattices and K3 surfaces; see
\cite{Bo} and its references.  $N_\WPi(\WJ)$ falls into 3 pieces: $\WJ$
itself, another Coxeter group $\WOmega$, and a group $\GOmega$ of
diagram automorphisms of $\WOmega$.  The last two groups are called
the ``reflection'' and ``non-reflection'' parts of the normalizer.
Borcherds bounded the virtual cohomological dimension of $\GOmega$ by
$|J|$.  Our theorems~\ref{thm-Wmin-dimension-bound},
\ref{thm-Borcherds} and \ref{thm-Wmax} give stronger bounds, in terms
of the types of the components of $J$ rather than the number of
nodes.  There are choices involved in the definition of $\WOmega$ and
$\GOmega$, and our bound in theorem~\ref{thm-Borcherds} applies
regardless of how these choices are made
(theorem~\ref{thm-Wmin-dimension-bound} is a special case).
Theorem~\ref{thm-Wmax} improves this bound when $\WOmega$ is
``maximal''.  In this case, when $J=D_5$ or $A_\modd$, $\GOmega$ turns
out to either be free or have an index~$2$ subgroup that is free.
This extends Brink's result \cite{Br} that $\GOmega$ is free when $J=A_1$.

The author is grateful to the Clay Mathematics Institute, the Japan Society
for the Promotion of Science, and Kyoto University for their support
and hospitality.

\bigskip
We follow the notation of \cite{Bo}, and refer to \cite{Humphreys} for
general information about Coxeter groups.  Suppose $(\WPi,\Pi)$ is a
Coxeter system, which is to say that $\WPi$ is a Coxeter group and
$\Pi$ is a standard set of generators.  The Coxeter diagram is the
graph whose nodes are $\Pi$, with an edge between $s_i,s_j\in\Pi$
labeled by the order $m_{ij}$ of $s_is_j$, when $m_{ij}>2$.  $\WPi$
acts isometrically on a real inner product space $V_\Pi$ with basis
(the simple roots)  $\Pi$ and inner products defined in
terms of the $m_{ij}$.  The (open) Tits cone $K$ is an open convex
subset of $V_\Pi^*$ on which $\WPi$ acts properly discontinuously with
fundamental chamber $\CPi$.  (Our $\CPi$ and $K$ are ``missing'' the
faces corresponding to infinite parabolic subgroups of $\WPi$.)  The
standard generators act on $V_\Pi^*$ by reflections across the
hyperplanes containing the facets of $\CPi$, and they also act on
$V_\Pi$ by reflections.
For a root $\a$ (i.e., a $\WPi$-image of a simple root) we write $\aperp$
for $\a$'s mirror, meaning the fixed-point set in $K$ of the
reflection associated to $\a$.

Now let $J\sset \Pi$ be a spherical subdiagram, i.e., one
corresponding to a finite subgroup of $\WPi$, and let $\Wmin$ be the
group generated by the reflections in $\WPi$ that act trivially on
$V_J\sset V_\Pi$.  This is the ``reflection'' part of $N_{\WPi}(\WJ)$,
or rather the strictest possible interpretation of this idea.  It
corresponds to Borcherds' $\WOmega$ in the case that the groups he
calls $\GPi$ and $\GJ$ are trivial; see the discussion after
lemma~\ref{lem-Wmin-contains-reflection}.  Let $\Jperp:=\cap_{\a\in
  J}\a^\perp$, pick a component $\Cmino$ of the complement of
$\Wmin$'s mirrors in $\Jperp$, and define $\Cmin$ as its closure (in
$\Jperp$).  By definition, $\Wmin$ is a Coxeter group, and the general
theory of these groups shows that $\Cmin$ is a chamber for it.  The
``non-reflection'' part of $N_\WPi(\WJ)$ means the subgroup $\Gmin$ of
$\WPi$ preserving $J$ (regarded as a set of roots) and sending $\Cmin$
to itself.  The reason for the first condition is to discard the
trivial part of $N_{\WPi}(\WJ)$, namely $\WJ$ itself.  That is,
$\Wmin{:}\Gmin$ is a complement to $\WJ$ in $N_{\WPi}(\WJ)$.  We write
$\Gmin^\vee$ for the subgroup of $\Gmin$ acting trivially on $J$
(equivalently, on $V_J$).  The reason for passing to this
(finite-index) subgroup is that $\Gmin$ often contains torsion and
therefore has infinite cohomological dimension for boring reasons.

\begin{theorem}
\label{thm-Wmin-dimension-bound}
$\Gmin^\vee$ acts
freely on a contractible cell complex of dimension at most 
\begin{equation}
\label{eq-Wmin-dimension-bound}
\#A_1 + \#D_{m>4} + \#E_6 + \#I_2(5)
+2\bigl(\#A_{m>1} + \#D_4\bigr)
\end{equation}
where $\#X_m$ means the number of components of $J$ isomorphic to a
given Coxeter diagram $X_m$.  In particular, $\Gmin^\vee$'s cohomological
dimension is at most \eqref{eq-Wmin-dimension-bound}.
\end{theorem}

Borcherds' result \cite[thm.~4.1]{Bo} has $|J|$ in place of
\eqref{eq-Wmin-dimension-bound}, but treats a more general group
$\GOmega$, of which $\Gmin$ is a special case.  The more general
case follows from this one, in theorem~\ref{thm-Borcherds} below.

\begin{proof}
\begingroup
\def\Gminx{\Gamma_{{\rm\! min},x}}
First we
prove for $x\in\Cmino$ that its stabilizer $\Gminx^\vee$ is trivial.
The $\WPi$-stabilizer of $x$ is some $\WPi$-conjugate $W_x$ of a spherical
parabolic subgroup of $\WPi$.  So $W_x$ acts on $V_\Pi$ as a finite Coxeter
group.  It is well-known that any vector stabilizer in such an action
is generated by reflections, so the subgroup $W_{x,J}$ fixing $J$
pointwise is generated by reflections.  Observe that any reflection in
$W_{x,J}$ lies in $\Wmin$.  Since $x$ lies in the interior $\Cmino$ of $\Cmin$,
it is fixed by no reflection in $\Wmin$, so
there can be no reflection in $W_{x,J}$, so $W_{x,J}=1$.  It is easy
to see that $W_{x,J}$ contains $\Gminx^\vee$, so we have proven that
$\Gmin^\vee$ acts freely on $\Cmino$.

$\Cmino$ is contractible because it is convex, and it obviously admits
an equivariant deformation-retraction to its dual complex.  So it
suffices to show that the dual complex has dimension at most
\eqref{eq-Wmin-dimension-bound}.  Suppose $\phi\sset\Jperp$ is a face
of a chamber of $\WPi$, with codimension in $\Jperp$ larger than
\eqref{eq-Wmin-dimension-bound}; we must show
$\phi\cap\Cmino=\emptyset$.  For some $w\in \WPi$, $w\phi$ is a face
of $\CPi$ whose corresponding set of simple roots $I'\sset\Pi$
contains $J':=w(J)\iso J$.  By the codimension hypothesis on $\phi$,
$|I'|-|J'|$ is more than \eqref{eq-Wmin-dimension-bound}.  Applying
the lemma below to $J'$ and $I'$, we see that $\WIp$ contains a
reflection $r$ fixing $J'$ pointwise.  Since $r\in\WIp$, its mirror
contains $w\phi$.  So $w^{-1}rw$ is a reflection fixing $J$ pointwise
(so it lies in $\Wmin$), whose mirror contains $\phi$.  Since $\Cmino$
is a component of the complement of the mirrors of $\Wmin$, it is disjoint from
$\phi$, as desired.
\endgroup
\end{proof}

\begin{lemma}
\label{lem-Wmin-contains-reflection}
If $J$ lies in a spherical Coxeter diagram $I\sset\Pi$, whose cardinality
exceeds that of $J$ by more than \eqref{eq-Wmin-dimension-bound}, then $\WI$
contains a reflection fixing $J$ pointwise.
\end{lemma}

\begin{remark}
Equality in \eqref{eq-Wmin-dimension-bound} holds when
$I$ extends the $A_m$, $D_m$, $E_6$ and $I_2(5)$ components of $J$ by
$A_1\to A_2$, $A_{m>1}\to D_{m+2}$, $D_4\to E_6$, $D_{m>4}\to
D_{m+1}$, $E_6\to E_7$ and $I_2(5)\to H_3$.  One can check in this
case that the conclusion of the lemma fails.
\end{remark}

\begin{proof}
We may suppose $I=\Pi$, by discarding the rest of $\Pi$.  Working one
component at a time, it suffices to prove the lemma under the
additional hypothesis that $\Pi$ is connected.  We now consider the
various possibilities for $\Pi$, and suppose $\WPi$ contains no
reflections fixing $V_J$ pointwise.  That is, we assume $\Wmin=1$.  In each case
we will show that $|\Pi|-|J|$ is at most \eqref{eq-Wmin-dimension-bound}.

The $\Pi=A_n$ case is a model for the rest.  Suppose the component of
$J$ nearest one end of $\Pi$ has type $A_m$ and does not contain that
end.  Then it must be adjacent to that end (since $\Wmin=1$), so together with the end it forms an $A_{m+1}$.  We
conjugate by the long word in $W(A_{m+1})$, which exchanges the two
$A_m$ diagrams in $A_{m+1}$ and fixes the roots in the other components
of $J$.  The result is that we may suppose without loss that $J$
contains that end of $\Pi$.  Repeating the argument to move the other
components of $J$ toward that end, we may suppose that there is
exactly one node of $\Pi$ between any two consecutive components of $J$.  And the
other end of $\Pi$ is either in $J$ or adjacent to it.  It is now
clear that $|\Pi|-|J|$ is the number of components of $J$, or one less
than this.  Since every component of $J$ has type~$A$, $|\Pi|-|J|$ is
at most~\eqref{eq-Wmin-dimension-bound}.
This finishes the proof in the $\Pi=A_n$ case.

If $\Pi=B_n=C_n$ then we begin by shifting any type $A$ components of
$J$ as far as possible from the double bond.  If $J$ has no $B_m$ then
$J$ contains one end of the double bond, and we get $|\Pi|-|J|$ equal
to the number of components of $J$, all of which have type~$A$.  If $J$ has a $B_m$ then the node after it (if
there is one) must be adjacent to some type $A$ component of $J$.
This is because $W(B_{m+1})$ contains a reflection acting trivially on
$V_{B_{m}}$.  This is easy to see in the model of $W(B_{m+1})$ as the
isometry group of $\Z^{m+1}$.  It follows that $|\Pi|-|J|$ is
the number of components of $J$ of type $A$.

In the $\Pi=D_{n>3}$ case, one can use the shifting trick to reduce to
one of the cases
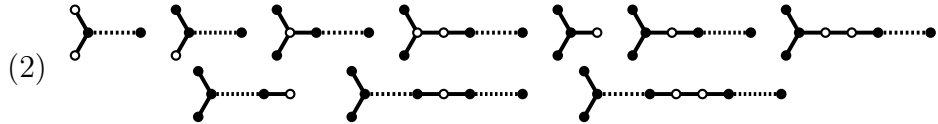
\psset{unit=\Dscaleunit}
\def\noderadius{\Dscalenoderadius}
\def\singlewidth{\Dscalesinglewidth}
\begingroup
\begin{equation}
\label{eq-Dn-diagrams}
\begin{matrix}
%
\begin{pspicture}(-58,-96)(208,96)
\single{-50}{-87}{0}{0}
\single{-50}{87}{0}{0}
\dashed{0}{0}{200}{0}
\hollow{-50}{-87}
\hollow{-50}{87}
\solid{0}{0}
\solid{200}{0}
\end{pspicture}
\quad
%
\begin{pspicture}(-58,-96)(208,96)
\single{-50}{-87}{0}{0}
\single{-50}{87}{0}{0}
\dashed{0}{0}{200}{0}
\hollow{-50}{-87}
\solid{-50}{87}
\solid{0}{0}
\solid{200}{0}
\end{pspicture}
\quad
%
\begin{pspicture}(-58,-96)(308,96)
\single{-50}{-87}{0}{0}
\single{-50}{87}{0}{0}
\single{0}{0}{100}{0}
\dashed{100}{0}{300}{0}
\solid{-50}{-87}
\solid{-50}{87}
\hollow{0}{0}
\solid{100}{0}
\solid{300}{0}
\end{pspicture}
\quad
%
\begin{pspicture}(-58,-96)(408,96)
\single{-50}{-87}{0}{0}
\single{-50}{87}{0}{0}
\single{0}{0}{200}{0}
\dashed{200}{0}{400}{0}
\solid{-50}{-87}
\solid{-50}{87}
\hollow{0}{0}
\hollow{100}{0}
\solid{200}{0}
\solid{400}{0}
\end{pspicture}
\quad
\begin{pspicture}(-58,-96)(108,96)
\single{-50}{-87}{0}{0}
\single{-50}{87}{0}{0}
\single{0}{0}{100}{0}
\solid{-50}{-87}
\solid{-50}{87}
\solid{0}{0}
\hollow{100}{0}
\end{pspicture}
\quad
\begin{pspicture}(-58,-96)(408,96)
\single{-50}{-87}{0}{0}
\single{-50}{87}{0}{0}
\single{0}{0}{200}{0}
\dashed{200}{0}{400}{0}
\solid{-50}{-87}
\solid{-50}{87}
\solid{0}{0}
\hollow{100}{0}
\solid{200}{0}
\solid{400}{0}
\end{pspicture}
\quad
\begin{pspicture}(-58,-96)(508,96)
\single{-50}{-87}{0}{0}
\single{-50}{87}{0}{0}
\single{0}{0}{300}{0}
\dashed{300}{0}{500}{0}
\solid{-50}{-87}
\solid{-50}{87}
\solid{0}{0}
\hollow{100}{0}
\hollow{200}{0}
\solid{300}{0}
\solid{500}{0}
\end{pspicture}
\\
\begin{pspicture}(-58,-96)(408,96)
\single{-50}{-87}{0}{0}
\single{-50}{87}{0}{0}
\dashed{0}{0}{200}{0}
\single{200}{0}{300}{0}
\solid{-50}{-87}
\solid{-50}{87}
\solid{0}{0}
\solid{200}{0}
\hollow{300}{0}
\end{pspicture}
\quad
\begin{pspicture}(-58,-96)(708,96)
\single{-50}{-87}{0}{0}
\single{-50}{87}{0}{0}
\dashed{0}{0}{200}{0}
\single{200}{0}{400}{0}
\dashed{400}{0}{600}{0}
\solid{-50}{-87}
\solid{-50}{87}
\solid{0}{0}
\solid{200}{0}
\hollow{300}{0}
\solid{400}{0}
\solid{600}{0}
\end{pspicture}
\quad
\begin{pspicture}(-58,-96)(808,96)
\single{-50}{-87}{0}{0}
\single{-50}{87}{0}{0}
\dashed{0}{0}{200}{0}
\single{200}{0}{500}{0}
\dashed{500}{0}{700}{0}
\solid{-50}{-87}
\solid{-50}{87}
\solid{0}{0}
\solid{200}{0}
\hollow{300}{0}
\hollow{400}{0}
\solid{500}{0}
\solid{700}{0}
\end{pspicture}
\end{matrix}
\end{equation}
\endgroup where the filled nodes are those in $J$ and the dashes
indicate a chain of nodes with no two
consecutive unfilled nodes.  (Except for the dashes on the left in the
last $3$ diagrams, which indicate chains of filled nodes.)  In every
case we get
$$
|\Pi|-|J|\leq \#A_1 + \#D_{m\geq4} + 2\,\#A_{m>1}.
$$
The most interesting case  is $A_{n-2}\to D_n$, at the top left.

We will treat the case $\Pi=E_8$ and leave the similar $E_6$ and $E_7$
cases to the reader.  If $J$ has a $D_4$, $D_5$ or $E_6$ component,
then it must also have a type $A$ component, and then $ |\Pi|-|J|\leq
2\,\#D_4 + \#D_5 + \#A_{m\geq1} $, as desired.  $J$ cannot be
$D_6$ or $E_7$, because then $\Wmin$ would contain
the reflection in the lowest root of $E_8$, which extends $E_8$ to the
affine diagram $\tildeE_8$.  So we may suppose $J$'s components have
type~$A$.  In order for $|\Pi|-|J|$ to exceed
\eqref{eq-Wmin-dimension-bound}, we must have $J=A_{m\leq5}$, $A_3A_1$,
$A_2A_1$ or $A_1^{m\leq3}$.  But none of these cases can occur,
because in each of them we may shift $J$'s components around so that
some node of $\Pi$ is not joined to $J$.

The remaining cases are $\Pi=F_4$, $H_3$, $H_4$ and $I_2$, the last
case including $G_2=I_2(6)$.  The facts required to treat these cases
are that if $J=B_2$ or $B_3$ in $\Pi=F_4$ then $\Wmin$ contains a
reflection, and similarly in the $J=H_3\sset H_4=\Pi$ case.  The first
fact is visible inside a $B_3$ or $B_4$ root system inside $F_4$.  To
see the second, observe that the root stabilizer in $H_4$ contains
Coxeter groups of types $A_2$ and $I_2(5)$, visible in the
centralizers of the two end reflections of $H_4$ (which are
conjugate).  So the root stabilizer can only be $W(H_3)$, which is to
say that the $H_3$ root system is orthogonal to a root.
\end{proof}

The greater generality obtained by Borcherds is the following. Let
$\GPi$ be a group of diagram automorphisms of $\Pi$, acting on $V_\Pi$
and $K$ in the obvious way.  The goal is to understand
$N_{\WPi{:}\GPi}(\WJ)$.  Again we discard the boring part of this
normalizer by passing to the subgroup $\WJ'$ preserving the set of
roots $J\sset\Pi$.  Let $\WOmega$ be any subgroup of $\WJ'$ which
contains $\Wmin$ and is generated by elements which act on $\Jperp$ by
reflections.  We define $\COmega^\circ$, $\COmega$ and $\GOmega$ as
for $\Cmino$, $\Cmin$ and $\Gmin$, and define $\GOmega^\vee$ as the
subgroup of $\GOmega\cap\WPi$ acting trivially on $J$.  (Borcherds
left $\GOmega^\vee$ unnamed and defined $\WOmega$ in terms of
auxiliary groups $R\trianglelefteq\GJ\sset\aut J$; his $\WOmega$ has
the properties assumed here.)  The inclusion $\Wmin\sset \WOmega$ is
the source of the subscript ``min'', but note that $\Cmin$ and $\Gmin$
are {\it larger\/} than $\COmega$ and $\GOmega$.  We can now recover
Borcherds' result \cite[thm.~4.1]{Bo} with our
\eqref{eq-Wmin-dimension-bound} in place of $|J|$.

\begin{theorem}
\label{thm-Borcherds}
Theorem~\ref{thm-Wmin-dimension-bound} holds with $\Gmin^\vee$ replaced by $\GOmega^\vee$.
\end{theorem}

\begin{proof}
The freeness of the action follows from the same argument.  (This is
why  $\GOmega^\vee$ is defined as a subgroup of $\GOmega\cap\WPi$
rather than just $\GOmega$.)  The essential point for the rest of
the proof is that $\WOmega$ contains $\Wmin$, so the decomposition of
$\Jperp$ into chambers of $\WOmega$ refines that of $\Wmin$.  This
shows $\COmega^\circ\sset\Cmino$.  So the dual complex of
$\COmega^\circ$ has dimension at most that of $\Cmino$, and we can
apply theorem~\ref{thm-Wmin-dimension-bound}.
\end{proof}

The point of considering $\WOmega$ rather than $\Wmin$ is that it is
larger and so $\GOmega$ will be smaller than $\Gmin$.  This is
good since the nonreflection part is more mysterious than the
reflection part.  So it is natural to define $\Wmax$ by setting
$\GPi=1$ and taking $\WOmega$ as large as possible, i.e., $\Wmax$
is the subgroup of $\WJ'$ generated by the transformations which act
on $\Jperp$ by reflections.  

This is the largest possible ``universal'' $\WOmega$, although a
larger $\WOmega$ is possible if $\Pi$ admits suitable diagram
automorphisms.  For example, $\GPi$ might contain elements acting on
$\CPi$ by reflections.  I don't know other examples, although
probably there are some.

We define $\Cmaxo$, $\Cmax$, $\Gmax$ and $\Gmax^\vee$ as above.  The
next theorem follows from lemma~\ref{lem-Wmax-contains-reflection} in
exactly the same way that theorem~\ref{thm-Wmin-dimension-bound}
follows from lemma~\ref{lem-Wmin-contains-reflection}.

\begin{theorem}
\label{thm-Wmax}
The dimension of the dual complex of $\Cmaxo$, hence the cohomological
dimension of $\Gmax^\vee$, is bounded above by 
\begin{equation}
\label{eq-Wmax-dimension-bound}
\#D_5 + 
\#A_{\modd} +
2\,\#A_{\meven}.
\end{equation}
\qed
\end{theorem}

\begin{remarks}
{\bf (i)} If $J$ has no $A_m$ or $D_5$ component then $\Gmax^\vee=1$ and
  $\Gmax$ is finite.  This is Borcherds' \cite[example~5.6]{Bo}.
{\bf (ii)} If $J=D_5$ or $A_\modd$ then $\Gmax^\vee\sset
  N_\WPi(\WJ)$ is free.  Also, since $|\aut J\,|\leq2$, $\Gmax^\vee$
  has index~$1$ or~$2$ in $\Gmax$. Therefore the non-reflection part
  $\Gmax$ of $N_\WPi(\WJ)$ has a free subgroup of index~$1$ or~$2$.
{\bf (iii)} If $J=A_1$ then $\Gmin=\Gmin^\vee=\Gmax=\Gmax^\vee$ has
cohomological dimension~${}\leq1$.  This recovers Brink's result
\cite{Br} that $\Gmin$ is free.
{\bf (iv)} If $J=A_\meven$ then the conclusion $\dim({\rm
  dual\ of\ }\Cmino)\leq2$ suggests that $\Gmax$ is often
comprehensible, like the $J=A_6$ example of \cite[example~5.4]{Bo}.
\end{remarks}

\begin{lemma}
\label{lem-Wmax-contains-reflection}
If $J$ lies in a spherical Coxeter diagram $I\sset\Pi$, whose
cardinality exceeds that of $J$ by more than
\eqref{eq-Wmax-dimension-bound}, then $\WI$ contains an element
preserving the set $J$ of roots and acting on $\Jperp$ by a
reflection.
\end{lemma}

\begin{proof}
This is essentially the same as for lemma~\ref{lem-Wmin-contains-reflection}, using the following
additional ingredients.  For example, when $I=D_n$ one can use them to
show that the 5th, 7th,
8th and 10th cases of \eqref{eq-Dn-diagrams} are impossible, while the first can
only occur when $n$ is even.

First, if $J=E_6\sset E_7=I$ then $\WI$
contains the negation of $V_I$, which we follow by the long word in
$\WJ$ to send $-J$ back to $J$.  The composition is the claimed
element of $\WI$.  The same argument applies if $J=I_2(5)\sset
H_3=I$.  

Second, if $J=A_{\modd}\sset D_{m+2}=I$ as in the first diagram of
\eqref{eq-Dn-diagrams}, then consider the long word in $\WI$.  It
negates $J$ and exchanges and negates the two simple roots in $I-J$.
Following this by the long word in $\WJ$ yields the claimed element
of $\WI$.  (When $m$ is even, the long word in $\WI$ negates the
simple roots in $I-J$ without exchanging them, so it doesn't act on
$\Jperp$ by a reflection.)

Third, if $J=D_{m\geq3}\sset D_{m+1}=I$ then consider the model of
$\WI$ as the group generated by permutations and evenly many negations
of $m+1$ coordinates, with $\WJ$ the corresponding subgroup for the
first $m$ coordinates.  Letting $\sigma$ be the negation of the last
two coordinates, and following it by the element of $\WJ$ sending
$\sigma(J)$ back to $J$, gives the claimed element of $\WI$.
\end{proof}

There is a nice geometrical interpretation of the freeness of $\Gmin$
in the case $J=A_1$, developed further in \cite{Allcock-centralizers}.
Namely, the natural map $\Cmino\to\Cmino/\Gmin\sset K/\WPi=\CPi$ is
the universal cover of its image.  The image is got by discarding all
the codimension~$2$ faces of $\CPi$ corresponding to even bonds in
$\Pi$, discarding all codimension~$3$ faces, and taking the component
corresponding to $J$.  This identifies $\Gmin$ with the fundamental
group of $J$'s component of the ``odd'' subgraph of $\Pi$ in a natural
manner.

One can extend this picture to the case $J\neq A_1$, but complications
arise.  First, one must take $\WOmega$ to be normal
in $\WPi{:}\GPi$.  Second, while
$\COmega^\circ\to\COmega^\circ/\GOmega^\vee$ is a covering space, the
image $\COmega^\circ/\GOmega$ of $\COmega^\circ$ in $\CPi$ is the
quotient of $\COmega^\circ/\GOmega^\vee$ by the finite group
$\GOmega/\GOmega^\vee$.  Usually,
$\COmega^\circ\to\COmega^\circ/\GOmega$ is only an orbifold cover
since $\GOmega$ often has torsion.  The top-dimensional strata of
$\COmega^\circ/\GOmega^\vee$ correspond to the ``associates'' of the
inclusion $J\to\Pi$ in the sense of \cite{Bo} and \cite{BH}.  Suppose
$J'\sset\Pi$ is (the image of) an associate and $I'$ is a spherical
diagram containing it.  Then the face of $\CPi$ corresponding to $I'$,
minus lower-dimensional faces, lies in $\COmega^\circ/\GOmega$
just if $\WIp$ contains no element preserving $J'$, acting on it in a
manner constrained by the choice of $\WOmega$, and acting on $\Jpperp$
by a reflection.  From this perspective,
lemmas~\ref{lem-Wmin-contains-reflection}
and~\ref{lem-Wmax-contains-reflection} amount to working out two cases
of Borcherds' notion of ``$R$-reflectivity''.  The orbifold structure
on $\COmega^\circ/\GOmega$ is essentially the same information as Borcherds'
classifying category for $\GOmega$.

\end{document}